\documentclass[10pt]{article}
 
\usepackage[latin1]{inputenc}
\usepackage[title]{appendix}
\usepackage{url, amsmath,enumerate,fancyhdr,amssymb, amsthm, 
pstricks, epsf, lscape, ifpdf, graphicx, float}
\usepackage{paralist}
\usepackage[margin=3.0cm]{geometry} 

\usepackage{algorithm}
\usepackage{algpseudocode}
\algnewcommand\algorithmicinput{\textbf{Input:}}
\algnewcommand\INPUT{\item[\algorithmicinput]}
\algnewcommand\algorithmicinitialization{\textbf{Initialization:}}
\algnewcommand\INITIALIZATION{\item[\algorithmicinitialization]}

\renewcommand{\arraystretch}{1.2}

\mathchardef\mhyphen="2D

\newtheorem{Definition}{Definition}

\newtheorem{Example}{Example}
\newtheorem{Proposition}{Proposition}
\newtheorem{Lemma}{Lemma}
\newtheorem{Theorem}{Theorem}
\newtheorem{Corollary}{Corollary}
\newtheorem{Remark}{Remark}
\newtheorem{Assumption}{Assumption}

\newcommand{\err}{\operatorname{err}}

\newcommand{\bpx}{\begin{pmatrix}}
\newcommand{\epx}{\end{pmatrix}}
\newcommand{\bbx}{\begin{bmatrix}}
\newcommand{\ebx}{\end{bmatrix}}

\newcommand{\bdef}{\begin{Definition}} 
\newcommand{\commentout}[1]{}
\newcommand{\co}[1]{}

\newcommand{\eref}[1]{(\ref{#1})}

\newcommand{\beq}{\begin{equation}}
\newcommand{\eeq}{\end{equation}}
\newcommand{\beqa}{\begin{eqnarray}}
\newcommand{\eeqa}{\end{eqnarray}}
\newcommand{\ba}{\begin{array}}
\newcommand{\ena}{\end{array}}
\newcommand{\bac}{\begin{array}{ccccccccccc}}
\newcommand{\eac}{\end{array}}
\newcommand{\bprop}{\begin{Proposition}}
\newcommand{\eprop}{\end{Proposition}}

\newcommand{\beqast}{\begin{eqnarray*}}
\newcommand{\eeqast}{\end{eqnarray*}}
\newcommand{\benum}{\begin{enumerate}}
\newcommand{\eenum}{\end{enumerate}}
\newcommand{\bit}{\begin{itemize}}
\newcommand{\eit}{\end{itemize}}
\newcommand{\bth}{\begin{Theorem}}
\newcommand{\enth}{\end{Theorem}}
\newcommand{\ble}{\begin{Lemma}}
\newcommand{\ele}{\end{Lemma}}
\newcommand{\bex}{\begin{Example}}
\newcommand{\eex}{\end{Example}}
\newcommand{\bcor}{\begin{Corollary}}
\newcommand{\ecor}{\end{Corollary}}
\newcommand{\brem}{\begin{Remark}}
\newcommand{\erem}{\end{Remark}}
\newcommand{\bass}{\begin{Assumption}}
\newcommand{\eass}{\end{Assumption}}

\renewcommand{\arraystretch}{1.2}

\newcommand{\bsmx}{\begin{small} \begin{pmatrix}}
\newcommand{\esmx}{\end{pmatrix} \end{small}}

\title{\Large  A simple preprocessing algorithm for semidefinite programming}

\author{Preston Faulk  \thanks{pfaulk@live.unc.edu, Department of Statistics and Operations Research, 
		University of North Carolina at Chapel Hill}  \hspace{1cm} G\'{a}bor Pataki  \thanks{gabor@unc.edu, Department of Statistics and Operations Research, 
		University of North Carolina at Chapel Hill }  \hspace{1cm} Quoc Tran-Dinh   \thanks{quoctd@email.unc.edu, Department of Statistics and Operations Research, 
		University of North Carolina at Chapel Hill } }

\begin{document}

\co{Change vs april 10: add section on generation}

\newcounter{algorithmcounter}

\maketitle 

\begin{abstract} 
	We propose a very simple preprocessing algorithm for semidefinite programming. 
	Our algorithm inspects the constraints of the problem, 
	deletes redundant rows and columns in the constraints, 
	and reduces the size of the variable matrix. It often detects infeasibility.
	Our algorithm does not rely on any optimization solver: the only subroutine it needs is Cholesky factorization, 
	hence it can be implemented with a few lines of code in machine precision. 
	We present computational results on a set of problems arising mostly 
	from polynomial optimization.
	\co{ and compare our method with the preprocessing algorithm of Permenter and Parrilo \cite{PerPar:14}, which relies on solving linear programming subproblems. }

\end{abstract} 

\section{Introduction and the preprocessing algorithm}\label{sec:intro}

Preprocessing is a key component in optimization solvers, in particular, in solvers of semidefinite programs (SDPs).  
SDPs -- whether they are formulated directly by a user, or whether they are the output of an algebraic modeling language -- often have  redundant constraints, or they may even be infeasible.
It is, of course, useful to detect these anomalies  
in a preprocessing stage.  
This paper provides a very simple preprocessing algorithm for SDPs, which can be implemented in a few line of codes in machine precision.

We consider  semidefinite programming problems in the form 

\begin{equation} \label{p} \tag{P}
\begin{array}{rrcl}
\inf & C \bullet X \\
s.t. & A_{i}\bullet X & = & b_{i} \,\,(i=1,\ldots,m) \\
& X & \succeq & 0
\end{array}
\end{equation}
where the $A_i$ are $n$ by $n$ 
symmetric matrices, the $b_i$ scalars, $X \succeq 0$ means that $X$ is symmetric, 
positive semidefinite (psd), 
and the $\bullet$ dot product of symmetric matrices is 
the trace of their regular product. We also write $X \succ 0$ to denote that $X$ is symmetric and positive definite. 

Our preprocessing algorithm is 
a very simple variant of facial reduction algorithms 
\cite{BorWolk:81, WakiMura:12, Pataki:13, LiuPataki:15}, which is also able to detect infeasibility. 
However, instead of solving SDP subproblems, like the algorithms proposed in these papers,  it  reduces the size of (\ref{p}) by simply inspecting the constraints.
\co{ like 
the methods in \cite{Friberg:16} and  \cite{GruberRendl:02}. }
Another related work is by Permenter and Parrilo in \cite{PerPar:14}, which finds reductions in SDPs by solving linear programming subproblems; our method is quite different though, since we do not rely on any optimization solver.
\co{
When (\ref{p}) is not {\em strictly feasible,} i.e., there is no feasible $X \succ 0, \,$ one can restate (\ref{p}) 
over a face of the semidefinite cone: this is the idea behind facial reduction algorithms 
\cite{BorWolk:81, WakiMura:12, Pataki:13}. In fact, facial reduction can even be used to detect infeasibility
\cite{WakiMura:12}. 
}

To motivate our algorithm, let us consider the following 
example:
\bex{\rm \label{ex1} 
The SDP instance (with an arbitrary objective function)
\begin{equation}  \label{mot-ex}
\begin{array}{cclcc}
\begin{pmatrix} 1 & 0 & 0 \\
0 & 0 & 0 \\
0 & 0 & 0 
\end{pmatrix} &\bullet&  X & = & 0 \\
\vspace{0.2cm}
\begin{pmatrix} 0 & 0 & 1 \\
0 & 1 & 0 \\
1 & 0 & 0 
\end{pmatrix} \, &\bullet& \, X & = & -1 \\
&& X & \succeq & 0,
\end{array}
\end{equation}
 is infeasible. Indeed, if $X = (x_{ij})_{i,j=1}^3$ were feasible in it, 
then $x_{11}=0,$   hence the first row and column of $X$ are zero
by positive semidefiniteness, so  the second constraint implies $x_{22} = -1, \,$ which is a contradiction.
}
\eex
Our preprocessing algorithm repeats the following basic step:

{\bf Basic step}
\bit
\item Find a constraint in (\ref{p}) which, after permuting rows and columns, and possibly multiplying both sides by $-1, \,$ 
is of the form 
 \beq \label{basic} 
 \bpx D & 0 \\ 0 & 0 \epx \bullet X = b, 
 \eeq
where $D \succ 0, b \leq 0.$ 
\item If $b < 0$ STOP; \eref{p} is infeasible.
\item If $b = 0, \,$ delete this constraint; also delete 
all rows and columns in the other constraints that correspond to rows and columns of $D.$ 
\eit

We stress that to find a constraint of the form \eref{basic}, 
we are only allowed to permute rows and columns in a constraint and to multiply 
both sides of a constraint by $-1;$ we do not take linear combinations of the constraints.

In Example \ref{ex1} in the first execution of the Basic step 
we find the first constraint, delete it, and also delete the first row and column from the second constraint matrix. In the second 
execution of the Basic step we find the constraint
$$
\bpx 1 & 0 \\ 0 & 0 
\epx \bullet X = -1,
$$
and declare that (\ref{mot-ex}) is infeasible. 

The idea of reducing SDPs by simply inspecting 
constraints already appears in several papers. Gruber and Rendl \cite{GruberRendl:02}, and Friberg \cite{Friberg:16} note that if $A \bullet X = 0$ is a constraint in an SDP with $A \succeq 0, \,$ then we can constrain $X$ to lie in a face of the positive semidefinite cone. These papers, however, do not assume 
the particular form of $A$ that we use and  they do not detect infeasibility.

Also, to the best of our knowledge, no such method has been implemented.

\section{Computational results}

We implemented our algorithm in Matlab, using 
incomplete Cholesky factorization 
to check positive definiteness. Since this subroutine of Matlab works in machine precision, which is of the order $10^{-16}$, our entire algorithm is implemented in machine precision.

Why should we preprocess SDPs? 
On the one hand, as our computational results show,
we cannot significantly reduce solution times on the tested instances.
On the other hand, we often detect their infeasibility.
Even when we do not, we bring the preprocessed instances closer to being strictly feasible, and strictly feasible problems behave better both from the theoretical, and the numerical point of view.
Since we do all computations with machine precision, we believe that 
we should always preprocess a problem, when we can.

Also, 
the constraints that we detect induce a fairly simple redundancy.
This redundancy should be clear even to a user not trained in 
optimization. Thus, returning the constraints that we find 
is  likely to help  him/her to better formulate other problems. 

We compare our preprocessing algorithm with the algorithm proposed by Permenter and Parrilo in \cite{PerPar:14}. Their algorithm solves linear programming subproblems to reduce the size of an SDP. It can work either on the problem \eref{p}, which we call the {\em primal}; or on its dual. 
It can use either diagonal, or diagonally dominant reductions -- for details, see
\cite{PerPar:14}. 

We calculated the DIMACS error measures (\cite{Mitt:13}) 
	for all instances before and after preprocessing. 
	In general, we will use the notation
	$$
	{\rm \err_{\rm before}} \, {\rm and} \, 	{\rm \err_{\rm after}}
	$$
	for the worst, i.e., largest DIMACS error measure before and after preprocessing
	for an instance. (Which particular instance we refer to will be clear from the context.)
	We also solved all instances using SDPT3 before and after preprocessing. 

We will say that the preprocessing (either ours, or the versions from \cite{PerPar:14}) {\em helped} a problem, if one of the following happens:

\begin{enumerate}
	\item It detects infeasibility, or
	\item $\err_{\rm before} > 10^{-6}$ and 
	\co{e worst DIMACS error measure before preprocessing is at least $10^{-6};$ and ratio }
	$$
\dfrac{\err_{\rm before}}{\err_{\rm after}} < \dfrac{1}{10};
	$$
	or 
	\item The optimal values before and after preprocessing differ by at least $10^{-6}.$
\end{enumerate}

We report our computational results in Table \ref{table1} on the problem set of 
$49$ instances from \cite{PerPar:14}. 
Our algorithm is denoted by FPT for the initials of the authors.
The algorithms in \cite{PerPar:14} are denoted by PP, and we also indicate 
whether they are applied to the primal or the dual problem, and whether they use diagonal (d) or diagonally dominant (dd) approximations.

In the first row we report the number of instances on which the algorithms found 
{\em some} reduction. In the second row we report the preprocessing time.
In the third  row we report the number of problems on which infeasibility is found, and fourth row the number of problems which are helped by the preprocessing.

It is important to note that the algorithms of Permenter and Parrilo do not detect infeasibility by themselves; to make a fair comparison we report if SDPT3 finds a problem infeasible after it has been preprocessed by one of their algorithms.

	\renewcommand*{\arraystretch}{1.3}
		\begin{table}[H]
		\begin{tabular}{l*{5}c}
					          & FPT &PP-Primal(d)& PP-Primal(dd) & PP-Dual(d)& PP-Dual(dd)  \\
					\hline
					\# Problems Preprocessed & 31 & 28 & 35 & 10 & 13  \\
					Total Preprocessing Time            & 33.42 & 30.630 & 96.02 & 6.20 &  89.29  \\
					Infeasibility Detect          & 14 & 15 ** & 16 ** & 1 ** & 0  \\
					Number Problems Helped     & 24 & 22 & 25   & 4  &  5   \\
				\end{tabular}
				\caption{Computational results. ''**'' means that SDPT3 detects infeasibility} 
				\label{table1} 
\end{table}
			
			The total time that SDPT3 took in solving the problems {\em before} 
			preprocessing is 253.8 seconds; the time it took after preprocessing is
			199.2 seconds. Thus the preprocessing only moderately helps in 
			reducing the solution time. 
			
In conclusion, our algorithm is competitive with the algorithms of \cite{PerPar:14} in terms of finding reductions, and detecting infeasibility.
At the same time it is simpler and  it does not rely on an optimization solver; 
thus we expect it to be at least as accurate, or more accurate.

We are currently exploring finding redundancies in (\ref{p}) by detecting 
constraints of the form 
$$
A \bullet X = b,
$$
where $A \succeq 0, \, b \leq 0.$ 

We prepared this version of the paper to be ready for the ICCOPT 2016 conference. In the final version we will add more tables and more details about  our computational results. 

\bibliographystyle{plain}
\bibliography{mysdp}
\end{document}